\documentclass[11pt,a4paper]{article}
 \usepackage[bottom=2.5cm,top=2.5cm,left=3cm,right=2.5cm]{geometry}
\usepackage{graphicx}
\usepackage{epstopdf}
\usepackage{psfrag}
\usepackage{amsmath,amssymb}
\usepackage{color}
\usepackage{nicefrac}
\usepackage{appendix}
\usepackage{floatrow}

\begin{document}


\title{Phase reconstruction with iterated Hilbert transforms}

\date{}
\author{Erik Gengel \footnote{Institute for Physics and Astronomy, Universit of Potsdam, Karl-Liebknecht-Str. 24/25, 14476 Potsdam, Germany} \footnote{Friedrich-Ebert Stiftung, egengel@uni-potsdam.de} and Arkady Pikovsky \footnote{Institute for Physics and Astronomy, Universit of Potsdam, Karl-Liebknecht-Str. 24/25, 14476 Potsdam, Germany} \footnote{Higher School of Economics, Nizhny Novgorod,
Russia}\footnote{Department of Control Theory, Institute of Information Technologies, Mathematics and Mechanics, Lobachevsky University Nizhny Novgorod, Russia, pikovsky@uni-potsdam.de}}
\maketitle
\section{Introduction and overview}

This chapter deals with the art of phase reconstruction. We focus on Hilbert transforms, 
however, much of the introduced methodology is not bound to Hilbert transforms alone. 

In general, three approaches to signal analysis of oscillatory signals
can be identified. The first approach applies statistical 
methods to extract information from observations assuming no further 
model \cite{li2010multiscale, schlemmer2018spatiotemporal, balzter2015multi}. The second approach 
takes the theory of dynamical systems into account and analyses the signals in terms of the 
phase and the amplitude notions provided in this 
theory \cite{Rosenblum-Pikovsky-01,Palus-Stefanovska-03,Cimponeriu2004,Bahr_etal-08,Rosenblum2006,%
Kralemann2007,kralemann2008phase,rosenblum2017reconstructing, topccu2018disentangling, Kralemann_etal-13}. 
In an intermediate methodology, a phase and 
an amplitude are extracted from the data and then analysed in terms of statistical quantities. These methods may or 
may not take an underlying theory into account \cite{yu2018functional, peng2005improved, rings2016distinguishing, stankovski2016alterations, zappala2019uncovering, paluvs2014multiscale, jajcay2018synchronization}
Alternatively, one applies machine learning 
techniques to obtain equations of motion directly from 
observations \cite{brunton2019data, yeung2019learning}. 
 
Here we focus on signal analysis approaches suitable for oscillating systems. 
The basic assumption is that the signal originates from a 
dynamical oscillating system, interacting with other systems and/or with the 
environment, and the goal is to understand the dynamics. 
This task is especially important and challenging in life science, where a theoretic description 
of the oscillators is in many cases lacking, 
because the underlying mechanisms are not clear. On the other hand, measurements of the 
full phase space dynamics are impossible, or 
would destroy the system itself. The latter aspect introduces the common setting where 
measurements of the systems are passive, i.e., an 
observer collects data from the free running system and may only apply 
weak perturbations 
to prevent damage. 
For such passive observations,
we pursue here the approach inspired by the dynamical system theory: we try to extract 
the phases from the signals, with the aim to build models as close to theoretical descriptions as possible.

The ideas of the phase dynamics 
reconstruction has been widely used in physics, chemistry, biology, medicine 
and other areas \cite{rosenblum2003synchronization, Stankovski_etal-17,blaha2011reconstruction,kralemann2011reconstructing} to understand properties 
of oscillators and coupling between them. The reason for this, as we discuss below, 
is that the phase is sensitive to interactions and external perturbations. In particular, 
many studies apply Hilbert transforms to reconstruct the phase from 
data (for example see \cite{benitez2001use, holstein2001synchronization, zappala2019uncovering, tort2010measuring} and references therein). 
However, several fundamental 
issues in the process of phase reconstruction are unresolved, 
long standing and mostly omitted in 
the community. One issue deals with the role of the 
amplitudes~\cite{rosenblum2018inferring, letellier1998non}.
 And from the view point of pure signal processing: How to deal 
 with phase-amplitude mixing in 
 Hilbert transforms \cite{guevara2017neural, cohen1999ambiguity}? The 
 latter issue will be 
 discussed in particular here and we describe a solution by virtue of \textit{iterative 
 Hilbert transform embeddings} (IHTE) \cite{gengel2019phase}. 

First, we describe the theoretical concepts. We then discuss the art of phase 
reconstruction with a 
focus on IHTE. We illustrate this method by presenting results for a Stuart-Landau 
non-linear 
oscillator, including reconstruction of the \textit{infinitesimal phase response curve} 
(iPRC). Finally, we discuss difficulties of application in case of noisy oscillations.

\section{Nonlinear oscillators and phase reduction}\label{sec: phase reduction}

Here we briefly review the phase reduction of driven limit cycle oscillators, 
for more details see~\cite{pikovsky2001synchronization,pietras2019network}.
An autonomous oscillator is described by $N$ state variables $\mathbf{y}$ which evolve
according to a system of differential equations 
$\dot{\mathbf{y}} = \mathbf{f}(\mathbf{y})$. One assumes
that this system has a stable limit cycle $\textbf{y}_0(t)=\mathbf{y}_0(t+T)$ describing
periodic (period $T$) oscillations. In the basin of attraction of the cycle one can always introduce a phase
variable $\varphi$ which grows uniformly in time 
\begin{equation}
\dot{\varphi} = \omega = \frac{2\pi}{T}\; . \label{eq asymp phase}
\end{equation} 
On the limit cycle, only the phase varies, so that $\textbf{y}_0(\varphi)=\textbf{y}_0(\varphi+2\pi)$,
which means that the value of the phase uniquely determines the point on the limit cycle.

If the autonomous oscillator is perturbed, i.e.\ 
it is driven by a small external force  
$\dot{\mathbf{y}}=\mathbf{f}(\mathbf{y})+\varepsilon\mathbf{p}(\mathbf{y},t)$, then the system 
slightly (of order $\sim\varepsilon$)
deviates from the limit cycle, and additionally the phase does not grow uniformly, but obeys (in the first order
in $\varepsilon$) the equation
\begin{equation}
\dot \varphi = \omega + \varepsilon Q(\varphi,t) \; , 
\label{eq:pert}
\end{equation}
where $Q$ can be expressed via $\mathbf{f},\mathbf{p}$ (see~\cite{pietras2019network} for details). 
Equation~\eqref{eq:pert}
contains only the phase and not the amplitude, it can be viewed as a result of 
the \textit{phase reduction}. The dynamics of the phase according to \eqref{eq:pert} 
allows for studying different important effects of synchronization, etc.
In the case when the oscillator is forced by another one, the force $\mathbf{p}(\eta)$
can be viewed as a function of the phase $\eta(t)$ 
of this driving oscillator, so the function $Q(\varphi,\eta)$ becomes
the coupling function depending on two phases. In experimental situations it is 
quite common to perturb just one variable of the system. In that case, if
the forcing term is scalar and does not depend on the system variables,
one can factorize $Q(\varphi,t)=Z(\varphi)P(t)$ into the iPRC $Z(\varphi)$ 
and the (scalar) external driving $P(t)$ 
\cite{winfree2001geometry, brown2004phase}. 
 
\begin{figure}[h!tbp]
\centering 
 \psfrag{S110}[cc][cc][1.0][0]{\raisebox{0.cm}{$X_{2,3}(\theta_{0,10})$,$X_{2,3}(\varphi)$}}
 \psfrag{phase}[cc][cc][1.0][0]{\raisebox{0.cm}{$\theta_{1,10}$,$\varphi$}}
 \psfrag{time}[cc][cc][1.0][0]{\raisebox{0.cm}{$t$}}
 \psfrag{X}[cc][cc][1.0][0]{\raisebox{0.cm}{$X_{2,3}(t)$}}
 \psfrag{a}[cc][cc]{\raisebox{-0.0cm}{\hspace{0cm}(a)}}
 \psfrag{b}[cc][cc]{\raisebox{-0.0cm}{\hspace{0cm}(b)}} 
\includegraphics[trim= 0cm 0.2cm 0cm 0.0cm, width=0.95\columnwidth]{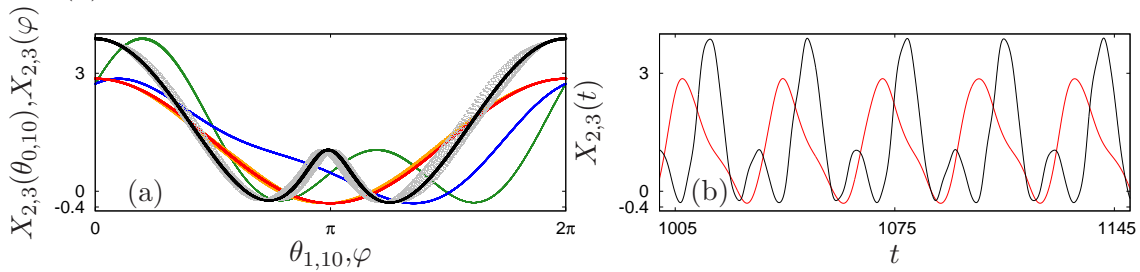}
\caption{Panel (a): Observables $X_2$ (blue points) and $X_3$ (green points) as functions of the true
phase $\varphi$. The fact that these sets
are not distiguishable from a line demonstrates validity of the phase description for the SL
oscillator (i.e. the amplitude modulation is indeed small).
The same observables as functions of the first protophase $\theta_1$ are not lines but broad sets (orange points for $X_2(\theta_{1})$ and grey points for  $X_3(\theta_{1})$). 
The same observables become good function of the protophase $\theta_{10}$ after 10-th iteration
of our procedure
(red and black points, corerspondingly). Panel (b): Time series for observables $X_{2,3}(t)$ (red,black). Simulation parameters are $\mu=8$, $\alpha=0.1$, $\nu=1$, $\varepsilon=0.1$ and $r=1.8$ (for $X_2(t)$) and $r=5.6$ (for $X_3(t)$). In this scale, small amplitude and phase modulations are hardly seen.
}
\label{fig 0.0}
\end{figure}

\paragraph{Example: Forced Stuart-Landau oscillator.}
In this contribution we consider as an example the perturbed Stuart-Landau oscillator (SL) 
\begin{equation}
\dot{a} = (\mu + i\nu)a - (1 + i\alpha)a|a|^2 + i\varepsilon P(t) \; ,\qquad P(t) = \cos(r\omega t)
\label{eq:slo}
\end{equation}
where $a(t) := R(t) \exp[i\phi(t)]$ is the complex amplitude. Parameter $\mu$ 
determines the amplitude ($\sqrt{\mu}$) and stability of the limit cycles, $\alpha$ 
is the nonisochronicity parameter. It is easy to check that
\begin{equation}
\varphi(t) = \phi(t) - \alpha \ln[R(t)]  \label{eq sl PHI}
\end{equation} 
is the proper phase, rotating, independently of amplitude $R$, with uniform frequency
 $\omega = \nu - \mu \alpha$. 
 The frequency of 
 the forcing is $r\omega$, where parameter $r$ is the ratio of the 
 external frequency to 
 the base frequency $\omega$.
 In the first order in $\varepsilon$, the amplitude and the phase dynamics read
 \begin{equation}
 \begin{aligned}
\dot{R} &= R(\mu-R^2) + \varepsilon P(t) \sin(\varphi)\;, \\
 \dot{\varphi} &= \omega +\varepsilon \mu^{-1/2}(\cos(\varphi) - \alpha \sin(\varphi))P(t)  \;.
 \end{aligned}
 \label{eq:slp}
 \end{equation}
Here, the iPRC is 
\begin{equation}
Z(\varphi) = (\cos(\varphi) - \alpha \sin(\varphi))\mu^{-1/2} \; .
\end{equation} 
One can see that for small $\varepsilon$ the dynamics of the SL is nearly 
periodic, with small ($\sim\varepsilon$) amplitude and phase modulations. Below
in this paper we will consider three different scalar observables of the SL dynamics:
$X_1(t)=\text{Re}[a(t)]$, $X_2(t) = 0.1(\text{Im}[a])^2+0.2(\text{Re}[a])^2+0.3\text{Im}[a]+0.4\text{Re}[a]$, and $X_3 = X_2 + 0.3\text{Re}[a]\text{Im}[a]$. The observable 
 $X_1$ is ``simple'',
 it is a pure cosine function of time for the autonomous SL oscillator. 
 The observable $X_2$ is also relatively simple (with one maximum and minimum pro period), 
 but not a pure cosine. The 
 observable $X_3$ can be viewed as a 
 \textit{multi-component signal}~\cite{feldman2011hilbert}, with two maxima and minima pro period.  Snapshots of the 
 time series of corresponding signals $X_{2,3}(t)$  are illustrated in Fig.\ \ref{fig 0.0}(b).

\section{Phase Reconstruction and iterative Hilbert transform embeddings}

\subsection{Waveform, phase and demodulation} \label{sec waveform phase modul}

In section \ref{sec: phase reduction} we introduced the phase dynamics 
concept for weakly perturbed oscillators.
It is based on the equations of the original oscillator's dynamics. 
In the context of \textit{data analysis}, one
faces a problem of the phase dynamics reconstruction solely
from the observations of a driven oscillator. From the time series of a scalar observable,
one wants to reconstruct the phase dynamics equation~\eqref{eq:pert}. 

The first assumption we make is that the phase modulation of the process observed is much 
stronger than the amplitude modulation. Although, according to the theory, amplitude 
perturbations appear already in the leading order $\sim\varepsilon$ 
(cf. Eq.~\eqref{eq:slp}), these variations could be small if the stability of the limit cycle 
is strong. Indeed, like example Eq.\ \eqref{eq:slp} shows, perturbations of 
the amplitude are inverse proportional to the stability of the limit cycle
$\sim \varepsilon/\mu$, and are additionally small for $\mu$ large. 
Thus, for the rest of this chapter we assume 
that the dynamical process under reconstruction is solely determined by the dynamics of the phase.

Generally, a time series emanates from an observable $X[\textbf{y}(t)]$ of 
the systems dynamics. According to the assumption above, we neglect amplitude 
modulation which means that we assume $\textbf{y}=\textbf{y}_0$,
so that the scalar signal observed is
purely phase modulated
\begin{equation}
X(t) = X[\textbf{y}_0(\varphi(t))] =: S(\varphi(t))\;.  \label{eq pm signal}
\end{equation}
Here a $2\pi$-periodic function $S(\varphi)=X[\textbf{y}_0(\varphi)]$ is unknown, 
we call it the waveform. 
The reconstruction problem for the signal $X(t)$ is that of
finding the waveform $S(\varphi)$ and the phase $\varphi(t)$. 
In Fig. \ref{fig 0.0}(a), we illustrate  these waveforms
for the observables  $X_{2,3}$ of the SL oscillator. Plotting $X_{2,3}$ as functions of $\varphi$ with dots,
one gets extremely narrow lines which indicate 
that for chosen large stability of the limit cycle the amplitude dynamics can be neglected 
and decomposition is possible. On the contrary, if the observed signals possess 
essential amplitude modulation, $X(\varphi)$ would look like a band. In that case the above
representation \eqref{eq pm signal} is not adequate.

We stress here that a
decomposition into the waveform and the phase is not unique. Indeed, let us 
introduce a new 
monotonous
``phase'' $\theta(t)$ 
according to an arbitrary transformation 
\begin{equation}
\theta=\Theta(\varphi),\quad \Theta(\varphi+2\pi)=\Theta(\varphi)+2\pi,\quad \Theta'>0.
\label{eq:ptoprtp}
\end{equation}
Then the signal can be represented as $X(t)=S(\Theta^{-1}(\theta))=\tilde{S}(\theta)$ 
with a new waveform $\tilde{S}=S \circ \Theta^{-1}$. 
Variables $\theta(t)$ are called protophases~\cite{Kralemann2007,kralemann2008phase}. 
Examples for mappings Eq.\ \eqref{eq:ptoprtp} are 
depicted in Fig.\ \ref{fig 0.1}. To see the difference between protophases
and true phase $\varphi(t)$, let us consider the non-driven, 
non-modulated dynamics. Here the phase $\varphi(t)$ 
grows uniformly $\dot{\varphi}=\omega$, while the 
protophase $\theta(t)$ grows non-uniformly, as
\begin{equation}
\dot{\theta}=\Theta'(\varphi)\omega=\omega\Theta'(\Theta^{-1}(\theta))=f(\theta)\;.
\label{eq:prtpd}
\end{equation}
However, having a protophase and the function $f(\theta)$ governing its dynamics, one can transform to the true
phase $\varphi(t)$ by inverting relation \eqref{eq:ptoprtp}: 
\begin{equation}
\frac{d\varphi}{d\theta}=\frac{1}{\Theta'(\varphi)}=\frac{\omega}{f(\theta)},\quad \varphi=\int_0^\theta \frac{\omega d\theta'}{f(\theta')}\;.
\label{eq:cl}
\end{equation}
Note that Eq.\ \eqref{eq:cl} is well defined as by 
construction, $\dot{\theta}=f(\theta)>0$.
In the case one observes driven oscillations, one approximately 
estimates $f(\theta)=\langle \dot{\theta}\rangle$,
see \cite{kralemann2008phase} for details.

According to the discussion above, one can perform the phase reconstruction of an 
observed signal $X(t)$ in two steps:\\
(i) Find a decomposition $X(t)=\tilde{S}(\theta(t))$ into a waveform and a 
protophase, satisfying conditions
\begin{equation}
\text{(I): }\; \forall t, \dot \theta > 0, \qquad \text{(II): }\; \tilde{S}(\theta) = \tilde{S}(\theta+ 2 \pi) \; .\label{eq constraints}
\end{equation}
(ii) Perform a transformation from a protophase to the phase, so that the latter grows on average uniformly
in time
\begin{equation}
\text{(III): } \langle \dot\varphi\rangle=\text{const}\;.
\label{eq:constr2}
\end{equation}
Conditions [I,II] ensure that the reconstructed protophase is monotonous and $2\pi$-periodic.  
Condition [III] selects the phase as a variable uniformly growing in time, in contrast to other 
protophases which according to \eqref{eq:prtpd} grow with a rate that is 
protophase-dependent (with $2\pi$-periodicity). Below we discuss in details 
the methods allowing for accomplishing steps (i) and (ii).

\subsection{Embeddings, Hilbert transform, and phase-amplitude mixing}
The first task, a decomposition into 
a waveform and a protophase, is trivial, if two scalar 
observables $\{X(t)=X[\mathbf{y}_0(t)],Y(t)=Y[\mathbf{y}_0(t)]\}$ of the 
oscillator's dynamics are available 
(of course, these observables should be not fully dependent). In this case, on the $\{X,Y\}$ plane one observes
a \textit{closed} continuous curve, parametrized by the phase, and the trajectory rotates along this curve. 
Any parametrization
of the curve, normalized by $2\pi$, will then provide a protophase as a function of time. After this, one has only
to accomplish the step (ii), i.e. to transform the protophase to the phase.

An intrinsically non-trivial problem appears, if only one scalar observable, $X(t)$, is available. The goal is
to perform a two-dimensional embedding of the signal $X(t)$, by generating from it the second variable $Y(t)$.
There exist several approaches for this task. The most popular ones are the 
delay-embedding $Y(t) = X(t-\tau)$ \cite{kim1999nonlinear}, the derivative 
embedding $Y(t) = \dot X(t)$ \cite{sauer1991embedology}, 
and the Hilbert transform (HT) embedding $Y(t) = \hat{H}[X](t)$, where (on a finite interval $[t_0,t_m]$)
\begin{equation}
\hat{H}[X](t) := \frac{\text{p.v.}}{\pi}\int_{t_0}^{t_m} \frac{X(\tau)}{t-\tau} d \tau\;. \label{eq hilbert time}
\end{equation} 

It is an observation of practice, that the latter approach based on the HT
often gives the most stable results. 
A reason for this is that the HT produces minimal distortions
to the signal's spectrum. Indeed, all the methods mentioned 
are linear transformations, which in Fourier space 
correspond to multiplications with 
factors $e^{i\Omega \tau}$, $i\Omega$, and $i\;\text{sign}(\Omega)$,
respectively. The factor for HT depends on frequency in a ``minimal'' way, and does not have, 
contrary to the delay embedding, a parameter. However, the HT provides only an approximate 
embedding, due to a mixing of phase and amplitude modulations~\cite{guevara2017neural}. 

Indeed, only for a non-modulated, i.e. for a purely periodic signal $X(t)$, the HT transform provides
a periodic $Y(t)$, so that on the $\{X,Y\}$ plane one observes a perfect closed loop. 
If the signal $X(t)$ is phase-modulated,
then on the $\{X,Y=\hat{H}[X(t)]\}$ plane one observes a non-closed trajectory 
(which only approximately can be considered as a loop), the width of the band gives the size of the appearing
amplitude modulation (see  Fig.~\ref{fig 0.0}(a) and Fig.~\ref{fig 0} below). 
(Also if one has a
purely amplitude-modulated signal, its HT will provide spurious phase modulation - but this is not
relevant for our problem). 
It should be noted that the 
spurious amplitude modulation arises solely due to the spectral properties of the 
Hilbert transform, and is not related to the length of the observation data. 
Usually, already 20 to 30 observed periods suffice to overcome 
boundary effects. Instead, the spectral content of the phase modulation 
heavily influences the appearance of amplitude
modulation, and hence the accuracy of reconstruction~\cite{gengel2019phase}.

In the next section we describe a method to circumvent this problem by 
virtue of \textit{iterated HT embeddings} (IHTE)~\cite{gengel2019phase},
illustrating the procedure with different observables of the SL oscillator.

\subsection{Iterated HT embeddings} \label{sec: IH}

\begin{figure}[h!tbp]
\centering 
 \psfrag{x1}[cc][cc][1.0][0]{\raisebox{0.cm}{$X_{1,2}(\theta_{0,10})$}}
 \psfrag{Hx1}[cc][cc][1.0][0]{\raisebox{0.cm}{$H[X_{1,2}](\theta_{0,10})$}}
 \psfrag{x2}[cc][cc][1.0][0]{\raisebox{0.cm}{$X_3(\theta_{0,10})$}}
 \psfrag{Hx2}[cc][cc][1.0][0]{\raisebox{0.cm}{$H[X_3](\theta_{0,10})$}}
 \psfrag{a}[cc][cc]{\raisebox{-0.0cm}{\hspace{0.0cm}(a)}}
 \psfrag{b}[cc][cc]{\raisebox{-0.0cm}{\hspace{0cm}(b)}} 
\includegraphics[trim= 0cm 0.2cm 0cm 0.5cm, width=0.95\columnwidth]{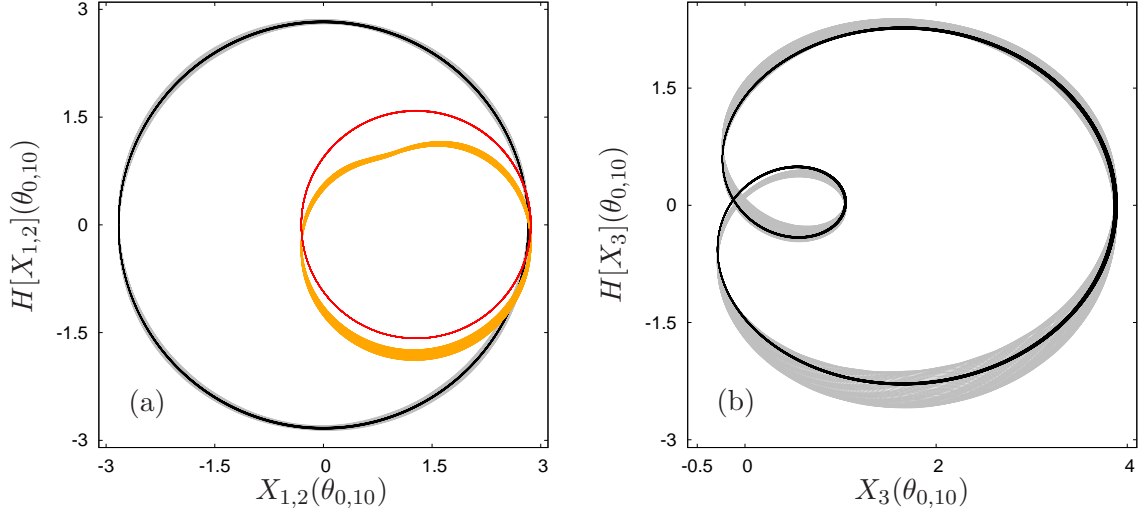}
\caption{IHTE for a periodically driven SL oscillator Eq.\ \eqref{eq:slo}
with harmonic driving $P(t)=\cos(r\omega t)$. 
Parameters: $\mu=8$, $\alpha=0.1$, $\nu=1$, $\varepsilon=0.1$. 
In panel (a) observables $X_1$ (for 
frequency ratio $r=5.6$) and $X_2$ (for $r=1.8$) are used. Shown are 
the first step of the IHTE hierarchy in grey and orange, and step ten in black and red 
for $X_1$ and $X_2$, respectively. In panel (b) the observable $X_3$ 
with $r=5.6$ is used where grey corresponds to the first embedding and black corresponds 
to the embedding in step ten.
Embeddings at the first iteration yield wide bands, which 
indicates for an ``artificial'' modulation of the amplitude, while at the 10th iteration the 
embeddings are nearly perfect lines, which means that the observed signals are nearly 
perfect phase modulated ones. Note that the embedding of $X_1$ 
has a circular shape, the embeddings of $X_{2,3}$ are distorted 
from a circle causing non-uniform protophases. In case of $X_3$, the embedding 
shows a loop (panel (b)).}
\label{fig 0}
\end{figure}

As discussed above, the HT embedding $\{X(t),\hat{H}[X(t)]\}$ although does 
not provide a closed looped line, allows one for an approximate determination of the 
protophase. To accomplish this, one needs to define a variable monotonously growing 
along the trajectory and gaining $2\pi$ at each approximate loop. A naive 
analytic-signal-based protophase $\text{arg}(X+iY)$ would work only for 
cosine-like waveforms like $X_{1,2}$. Therefore we employ another 
definition of the protophase, based on the trajectory length~\cite{kralemann2008phase} 
\begin{equation}
L(t) = \int_0^t \sqrt{\dot{X}^2(\tau)+\dot{Y}^2(\tau)} d\tau \label{eq: length } \; .
\end{equation}
This length grows monotonously also in the case when the
embedding has loops (cf. Fig.~\ref{fig 0}), in which case the
analytic-signal-based definition obviously fails.

Having calculated the length $L(t)$, we can transform it to a protophase by interpolation. 
For this, we define in the signal features which we attribute to the zero (modulo $2\pi$) 
protophase, and define the corresponding time instants $t_j$.
In the simplest case, one can define a protophase $\theta(t)$  on the interval $(t_j,t_{j+1})$ 
as a linear function of the length
$\theta(t)=2\pi j+2\pi (L(t)-L(t_j))/(L(t_{j+1})-L(t_j))$. However, such a protophase will be 
discontinuous in the first derivative. A better transformation is achieved via splines: 
one constructs a spline approximation for the function $\theta(L)$, provided one knows 
the values of this function at the signal features: $\theta(t_j)=2\pi j$ at $L(t_j)$. 

Constructed in this way, the protophase $\theta(t)$ is only approximate, 
because $X(\theta+2\pi)\neq X(\theta)$. Visually, on the plane $\{X,\hat{H}[X]\}$ one 
observes a band instead of a single loop (see Fig.~\ref{fig 0}).  Also, when $X$
is plotted vs $\theta$, one observes not a single-valued function, but 
a band (see Fig.~\ref{fig 0.0}(a)).

\begin{figure}
\includegraphics[width=0.7\textwidth]{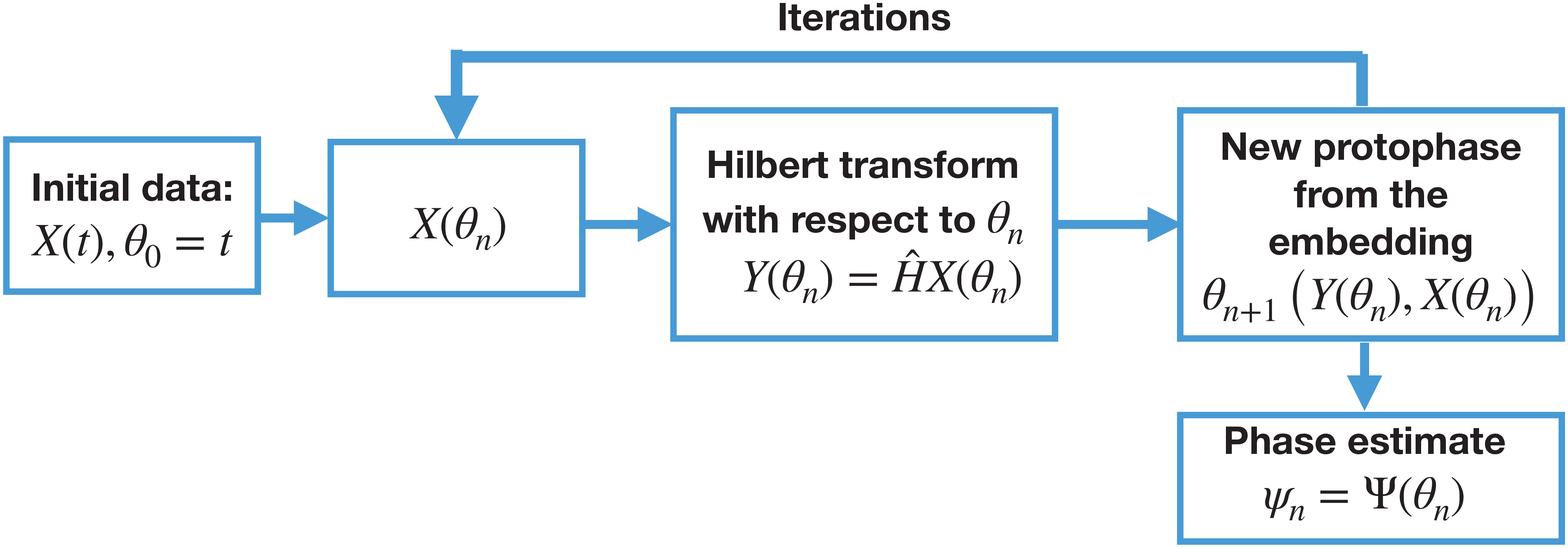}
\caption{
Here we schematically explain the 
iterative Hilbert tranform 
embeddings. Typically only one iteration is performed, and the 
protophase $\theta_1$ is used 
for further
analysis. We show in this chapter, how the quality of the phase reconstruction improves 
with the iterative embeddings.} 
\label{fig: 0.00}
\end{figure}

Recently, in Ref.~\cite{gengel2019phase}, we proposed to use iterative 
Hilbert transform embeddings (IHTE) to improve the quality of the protophase definition above. 
Our idea is to perform subsequent Hilbert transforms based on the previously 
calculated protophases $\theta_n(t)$, where $n$ denotes the step of iteration (see Fig.~\ref{fig: 0.00}).
Intuitively, the advantage of iterations can be understood as follows:
The widely used first iteration already presents an approximation to the protophase, although
not a perfect one. This means, that the function $X(\theta_1)$ still has modulation, but less than $X(t)$.
Now, if we take $\theta_1$ as a new time and again perform a demodulation by virtue of the
Hilbert transform embedding, we expect $\theta_2$ to be better than $\theta_1$, etc.
A detailed analysis performed in Ref.~\cite{gengel2019phase} shows that this procedure indeed converges to perfect demodulation.
 
In terms of iterations, the protophase $\theta(t)$ discussed above is the first 
iteration $\theta_1(t)$, while the time variable can be considered as the ``zero'' 
iteration $\theta_0(t)$. At each iteration step
we use the obtained protophase  as a new ``time'' with respect to which the next HT is performed:
\begin{equation}
Y_{n+1}(\theta_n)=\hat{H}[X(\theta_n)] := \frac{\text{p.v.}}{\pi}\int_{\theta_n(t_0)}^{\theta_n(t_m)} \frac{X(\theta'_n)}{\theta_n-\theta_n'} d \theta_n' \; . \label{eq: it hilbert}
\end{equation}

An implementation of this integral is given in \cite{gengel2019phase}. Basically there are two 
challenges here: first, the integration has to be performed on a non-uniform grid and second, 
one has to take care of the singularity at $\theta_n'=\theta_n$.
 
The iteration process will be as follows (see Fig.~\ref{fig: 0.00}):
\begin{enumerate}
\item[1.] Having $X(\theta_n)=X(t(\theta_n))$, we calculate 
$Y_{n+1}(\theta_n) = \hat{H}[X(\theta_n)]$ according to \eqref{eq: it hilbert}. 
\item[2.] Next, we construct the embedding $\{X,Y_{n+1}\}$ and find the 
length $L(\theta_n)$ from \eqref{eq: length }. 
\item[3.] After defining signal features, we calculate, using splines, 
the new protophase $\theta_{n+1}$ as a function of $L(\theta_n)$, which gives the 
new protophase $\theta_{n+1}$ as a function of the old one $\theta_n$.
\end{enumerate}

The steps 1-3 are repeated, starting from $\theta_0=t$. After $n$ iterations, 
we obtain a waveform and a protophase 
\begin{equation}
\tilde{S}(\theta_n) = X(t(\theta_n)) 
\end{equation}
As has been demonstrated in Ref.~\cite{gengel2019phase}, the procedure 
converges to a proper protophase, fulfilling conditions [I,II] above. 
For a purely phase modulated signal, at large $n$ the errors \eqref{eq:errs} 
reach very small values limited by 
accuracy of integration. The convergence rate depends heavily on the 
complexity of the waveform and on the level and frequency of modulation, 
but typically at $\hat{n}\approx 10$ a good protophase is constructed.

\begin{figure}[h!tbp]
\centering 
 \psfrag{prp}[cc][cc][1.0][0]{\raisebox{0.cm}{$\theta_{1,10}(\varphi)$,$\psi_{10}(\varphi)$}}
 \psfrag{vp}[cc][cc][1.0][0]{\raisebox{0.cm}{$\varphi$}}
 \psfrag{a}[cc][cc]{\raisebox{-0.0cm}{\hspace{3.0cm}(a)}}
 \psfrag{b}[cc][cc]{\raisebox{-0.0cm}{\hspace{3.0cm}(b)}} 
\includegraphics[trim= 0cm 0.5cm 0cm 0.5cm, width=0.95\columnwidth]{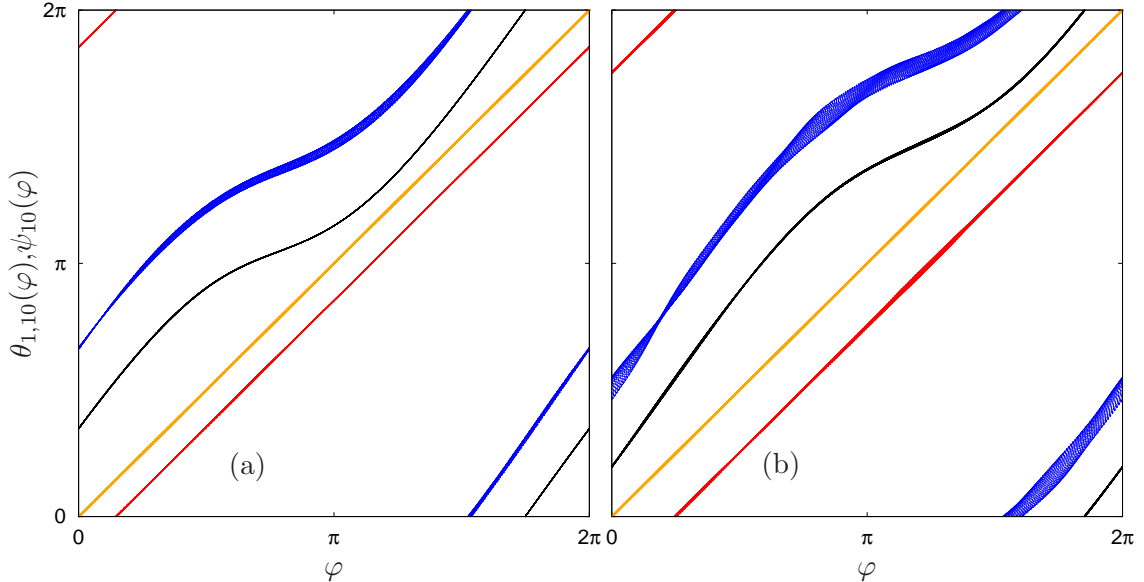}
\caption{
Depicted are the phase-to-protophase maps Eq.\ \eqref{eq:ptoprtp} for 
$X_2[.]$ (panel (a)) and $X_3[.]$ 
(panel (b)) based on the embeddings shown in Fig.\ \ref{fig 0}. Colours correspond 
to $\Theta_1(\varphi)$ (blue, this
data form a rather wide band indicating that the protophase at the first iteration is not precise), 
$\Theta_{10}(\varphi)$ (black,
this data forms a narrow line indicating for a good protophase reconstruction), $\psi_{10}(\varphi)$ (red,
this narrow line is straight indicating for a good phase reconstruction). The orange 
line is the diagonal. For better visibility  the curves are shifted vertically.}
\label{fig 0.1}
\end{figure}

Summarizing, the IHTE solve the problem of constructing a protophase 
$\theta(t)=\theta_{\hat{n}}(t)$ and the corresponding waveform $\tilde{S}(\theta)$ from a scalar 
phase-modulated signal $X(t)$; this protophase fulfills 
conditions \eqref{eq constraints}-[I,II]. Indeed, one observes in Fig.\ \ref{fig 0.1} that the 
first mapping $\Theta_1(\varphi)$ is not purely $2\pi$-periodic (blue bands). Instead, after 
ten iterations, $\Theta_{10}(\varphi)$ effectively has become a line (black) indicating that a 
protophase is reconstructed. The same can be seen in Fig.~\ref{fig 0.0}(a), where 
bands of values $X_{2,3}(\theta_1)$ are transformed to narrow lines $X_{2,3}(\theta_{10})$
after ten iterations.

As the final step in obtaining a close 
estimate $\psi(t)$ of the proper phase $\varphi(t)$, we have to perform 
the protophase-to-phase transformation, as described in Ref.~\cite{kralemann2008phase}. 
The transformation is based on relation~\eqref{eq:cl}, where the Fourier components of 
the density of the protophase are estimated according to 
$F_k=t_m^{-1}\int_0^{t_m}\exp[-ik\theta(t)]\;dt$; these components 
are used to perform the transformation as 
$\psi=\theta+\sum_{k\neq 0} F_k(ik)^{-1}[\exp(ik\theta)-1]$. Indeed, 
one observes in Fig.\ \ref{fig 0.1} (red lines) that $\psi(t)$ is, up to estimation errors, resembling the 
dynamics of $\varphi(t)$. However, we want to stress here that determination
of the protophase-to-phase transformation is based on a statistical evaluation
of the probability density of the protophase. Hence, in order to achieve a proper 
reconstructions with small distortions in the protophase-to-phase mapping, one 
needs long time series.

We can check for the similarity of $\theta_n(t)$ or $\psi_n(t)$ to the true phase $\varphi(t)$ by 
calculating a phase and a frequency error as the standard deviations
\begin{equation}
\begin{gathered}
\text{STD}^{\mathbf{q}}_n = \sqrt{\frac{1}{\hat{N}_1}\int_{t_{min}}^{t_{max}} [{\mathbf{q}}_n(\tau) - \varphi (\tau)]^2 d \tau} \qquad  \text{STD}^{\dot{{\mathbf{q}}}}_n = \sqrt{\frac{1}{\hat{N}_2}\int_{t_{min}}^{t_{max}} [\dot{{\mathbf{q}}}_n(\tau) - \dot{\varphi} (\tau)]^2 d \tau} \\
\hat{N}_1 = \int_{t_{min}}^{t_{max}} (\varphi(\tau)-\tilde{\omega} \tau)^2 d \tau \qquad \qquad \qquad \hat{N}_2 = \int_{t_{min}}^{t_{max}} [\dot{\varphi}(\tau)-\tilde{\omega}]^2 d \tau \; .
\end{gathered}
\label{eq:errs}
\end{equation}
STD$^{{\mathbf{q}},\dot{{\mathbf{q}}}}_n$ tend to zero only, if the reconstructed 
protophases and transformed protophases $\mathbf{q}_n=\{\theta_n(t),\psi_n(t)\}$ are 
close to the true phase $\varphi(t)$ of the system (see Eq.\eqref{eq sl PHI}). 
In the integration, we skip the outer ten percent at the beginning and at the end of the time series, 
to avoid boundary effects. Estimations of the instantaneous frequency $\dot{\varphi}$ and 
$\dot{\mathbf{q}}_n$ are performed by a 12th order polynomial filter (Savitzky-Golay filter) 
with a window of 25 points and four times repetition \cite{savitzky1964smoothing} denoted 
as SG[12,25,4].
Throughout the chapter we use a sampling rate of 
$dt=0.01$, such that the smoothing window has a width of $dt=0.25$, corresponding 
to roughly $11\%$ of the fastest forcing period ($r=14.3$).
The estimated average 
growth rate $\tilde{\omega}$ is obtained by linear regression. Note that the normalization integral $\hat{N}_1$ 
is suitable for all phases where the average growth is linear.

\section{Numerical experiments}

\subsection{Deterministic oscillations}\label{sec determ osci}

\begin{figure}[h!tbp!]
\centering 
 \psfrag{time}[cc][cc][1.0][0]{\raisebox{0.cm}{time}}
 \psfrag{ymod}[cc][cc][1.0][0]{\raisebox{0.cm}{$q(t)$,$u_{1,10}(t)$}}
 \psfrag{yIF}[cc][cc][1.0][0]{\raisebox{0.cm}{$\dot{\varphi}(t)$,$\dot{\theta}_{1,10}(t)$}} 
 \psfrag{yIFc}[cc][cc][1.0][0]{\raisebox{0.cm}{$\dot{\varphi}(t)$,$\dot{\psi}_{1,10}(t)$}} 
 \psfrag{d}[cc][cc]{\raisebox{-0.4cm}{\hspace{-0.2cm}(a)}} 
 \psfrag{e}[cc][cc]{\raisebox{-0.4cm}{\hspace{-0.2cm}(b)}} 
 \psfrag{i}[cc][cc]{\raisebox{-0.4cm}{\hspace{-0.2cm}(c)}} 
 \psfrag{j}[cc][cc]{\raisebox{-0.4cm}{\hspace{0.85cm}(d)}} 
\includegraphics[width=\columnwidth]{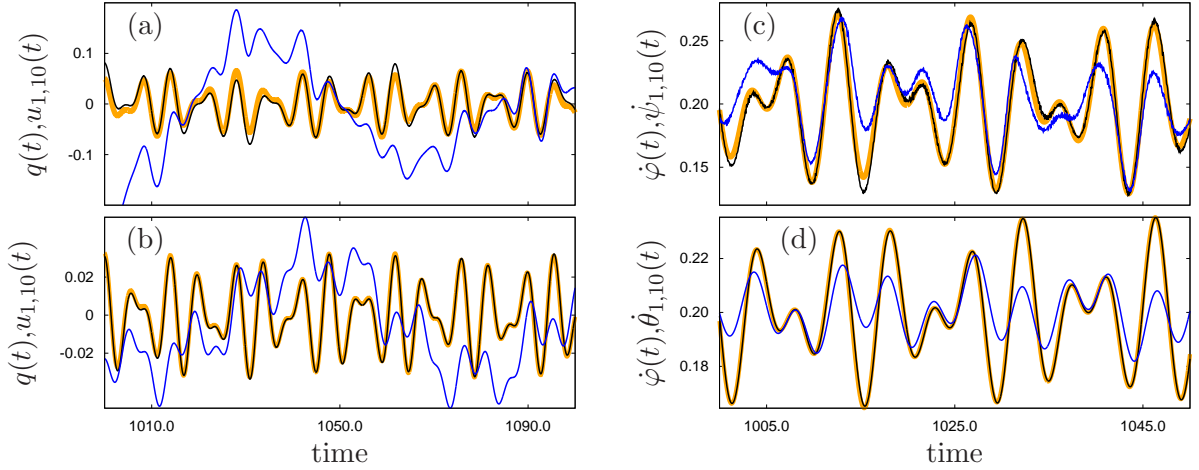}
\caption{
Comparison of true modulation $q(t)$ and of true instantaneous frequency
$\dot{\varphi}$ (orange) with the corresponding 
iteration results $u_n(t)$ in the first step (blue) and in the 10th step (black) for 
the observables $X_1$ (b,d) and $X_3$ (a,c). Parameters 
are $\mu=8$, $\alpha=0.1$, $\nu=1$, $\varepsilon=0.1$ and $r=5.6$. For 
calculation of the derivative $\dot{\varphi}$ we use a SG[12,25,4] filter.}
\label{fig 2}
\end{figure}

\begin{figure}[h!tbp!]
\centering 
 \psfrag{xe}[cc][cc][1.0][0]{\raisebox{0.cm}{$n$}}
 \psfrag{ye1}[cc][cc][1.0][0]{\raisebox{0.cm}{STD$^{\theta}_n$}}
 \psfrag{ye11}[cc][cc][1.0][0]{\raisebox{0.cm}{STD$^{\dot{\theta}}_n$}}
 \psfrag{ye2}[cc][cc][1.0][0]{\raisebox{0.cm}{STD$^{\psi}_n$}}
 \psfrag{ye21}[cc][cc][1.0][0]{\raisebox{0.cm}{STD$^{\dot{\psi}}_n$}}
 \psfrag{000c125}[cc][cc][0.75][0]{\raisebox{0.cm}{\hspace{0.18cm}$0.4$}}
 \psfrag{c127}[cc][cc][0.75][0]{\raisebox{-0.0cm}{\hspace{-0.41cm}$1.8$}}
 \psfrag{c128}[cc][cc][0.75][0]{\raisebox{-0.0cm}{\hspace{-0.41cm}$3.5$}}
 \psfrag{c129}[cc][cc][0.75][0]{\raisebox{-0.0cm}{\hspace{-0.41cm}$7.2$}}
 \psfrag{c130}[cc][cc][0.75][0]{\raisebox{-0.0cm}{\hspace{-0.6cm}$14.3$}}
 \psfrag{c20}[cc][cc][0.75][0]{\raisebox{-0.0cm}{\hspace{-0.7cm}$5.6$}}
 \psfrag{a}[cc][cc]{\raisebox{5.0cm}{\hspace{7.8cm}(a)}}
 \psfrag{b}[cc][cc]{\raisebox{5.0cm}{\hspace{7.8cm}(b)}} 
 \psfrag{c}[cc][cc]{\raisebox{5.0cm}{\hspace{7.8cm}(c)}} 
 \psfrag{d}[cc][cc]{\raisebox{5.0cm}{\hspace{7.8cm}(d)}} 
 \psfrag{e}[cc][cc][0.75][0]{\raisebox{15.3cm}{\hspace{13.3cm}$r=$}}
\includegraphics[width=\columnwidth]{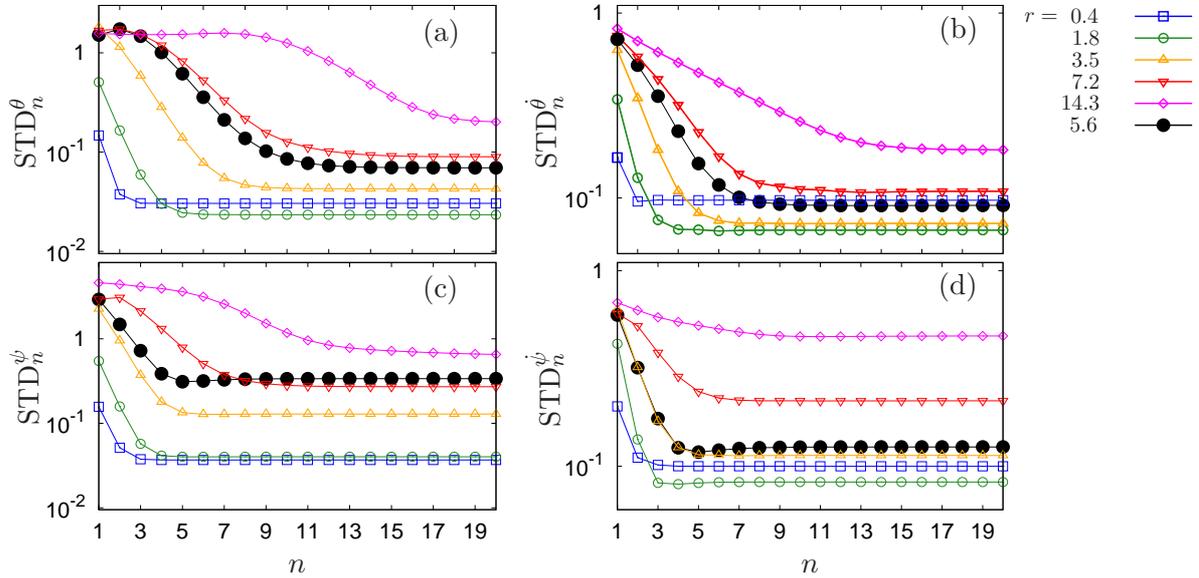}
\caption{
Phase and frequency errors for observables $X_1$ (a,b) 
and $X_2$ (c,d) for different forcing frequencies 
$r\omega$. Also shown in (c,d black doted line) is the reconstruction 
error for $X_3$, where the use of $L(t)$ for phase calculation is crucial. 
Slow modulations are essentially reconstructed in the first step while fast 
modulations need at least several iterations. 
With increased forcing frequency, $\theta_n(t)$ differs 
significantly from $\varphi(t)$ and the number of needed iterative steps grows.}
\label{fig 1}
\end{figure}

Here we consider the SL system~\eqref{eq:slo} with $\mu=8$, $\alpha=0.1$, $\nu=1$.
As the observables we explore functions $X_{1,2,3}[a]$ defined above. 
The system is forced harmonically 
by $\varepsilon P(\eta) = \varepsilon \cos(\eta(t))$ with amplitude $\varepsilon=0.1$. 
The external force phase is $\eta(t) = r \omega t$, for the explored range of driving frequencies
$r \omega$ the SL operated
in the asynchronous regime. We observe 100 periods with a time step of $dt=0.01$.

In Fig.\ \ref{fig 1} the phase and the frequency errors according to Eq.\ \eqref{eq:errs} 
for the first $20$ iteration steps are shown. While for slow modulations ($ r< 1$), the 
reconstruction is already accurate in the first step, for fast forcing frequencies ($r > 1$) 
indeed several iterations are needed for precise reconstruction. The reason for this is 
that for high-frequency modulations iterative HT embeddings first shift high-frequency 
Fourier components of the phase modulation to lower frequencies, where they eventually 
disappear. This mechanism is closely related to the Bedrosian identities \cite{xu2006bedrosian} 
and is explained in detail in \cite{gengel2019phase}. For the reconstruction of phases in 
case of $X_{2,3}[a]$, we have to calculate the transformed phase $\psi(t)$, 
because here the protophases deviate from uniform growth. The results show, that IHTE 
combined with the protophase-to-phase transformation provides proper phase reconstructions 
for the fairly stable limit cycle oscillator under study. 

Figure \ref{fig 2} presents comparisons of the inferred modulation 
$u_n(t):=\psi_n(t)- \tilde{\omega} t$  with
the true one $q(t)=\varphi-\omega t$,
and of the inferred instantaneous frequencies $\dot{\theta}_n(t)$/$\dot{\psi}_n(t)$ 
with $\dot\varphi$, for 
a quite fast external force $r=5.6$ (black bold dots in Fig.\ \ref{fig 1}). While the 
first iterate is by far not accurate, iterations provide the reconstructed estimation 
of the phases $\psi_{20}(t)$ which is very close to $\varphi(t)$.

\subsection{Reconstruction of the phase response curve from observation}

\begin{figure}[!htbp!] 
\floatbox[{\capbeside\thisfloatsetup{capbesideposition={right,top},capbesidewidth=5cm}}]{figure}[\FBwidth] 
{\caption{
The phase response curve $Z(\varphi)$ (orange). The blue unevenly dashed line 
depicts the estimation $\tilde{Z}_{30}(\varphi)$ based on $\theta_1(t)$. Also 
shown are the estimations $\tilde{Z}_{30}(\varphi)$ (solid line) 
and $\tilde{P}_{30}(\eta)$ (dashed line) based on $\theta_{10}(t)$ for $r=[0.06,4.5,5.6]$ 
(top to bottom). Red lines refer to the coupling functions Fig.\ \ref{fig 3}}\label{fig 3.1}
}
{ \psfrag{prcforce}[cc][cc][0.85][0]{\raisebox{0cm}{\hspace{0cm} $\tilde{Z}_{30}(\varphi),\tilde{P}_{30}(\eta)$}}
\psfrag{pprc}[cc][cc][0.8][0]{\raisebox{0.3cm}{$\varphi$}} 
\psfrag{e}[cc][cc][0.8]{\raisebox{0cm}{\hspace{0cm} (e)}} 
 \includegraphics[scale=1.2]{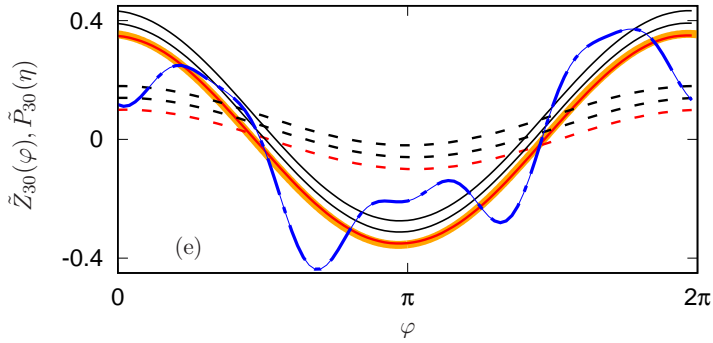}}
\end{figure}

\begin{figure}[h!tbp] 
\floatbox[{\capbeside\thisfloatsetup{capbesideposition={right,top},capbesidewidth=5cm}}]{figure}[\FBwidth] 
{\caption{
Reconstructed coupling functions $\tilde{Q}_{30,1}(\varphi,\eta)$ (panel (a)) 
and $\tilde{Q}_{30,10}(\varphi,\eta)$ (panel (b)) 
based on $\theta_1(t)$ and $\theta_{10}(t)$, 
respectively. The reconstruction error for the first iteration 
$Q_{30,1}(\varphi,\eta)-Q(\varphi,\eta)$ is shown in panel (c), 
and for the 10th iteration $Q_{30,10}(\varphi,\eta)-Q(\varphi,\eta)$ 
in panel (d). Noteworthy,
the vertical scale in panel (d) is more than ten times smaller than in panel (c).}\label{fig 3}
}
{ \psfrag{a}[cc][cc][0.8]{\raisebox{-1.5cm}{\hspace{0cm} (a)}} 
 \psfrag{b}[cc][cc][0.8]{\raisebox{-1.5cm}{\hspace{0cm} (b)}} 
 \psfrag{c}[cc][cc][0.8]{\raisebox{-1.5cm}{\hspace{0cm} (c)}} 
 \psfrag{d}[cc][cc][0.8]{\raisebox{-1.5cm}{\hspace{0cm} (d)}} 
 \psfrag{ap}[cc][cc][0.8][0]{\raisebox{1.3cm}{$\varphi$}} 
 \psfrag{ep}[cc][cc][0.8][0]{\hspace{-1.3cm}$\eta$} 
 \psfrag{PhPs}[cc][cc][0.85][0]{\raisebox{0.0cm}{$\varphi$,$\eta$}}
 \includegraphics[trim= 0cm 0.2cm 0cm -0.2cm, scale=0.9]{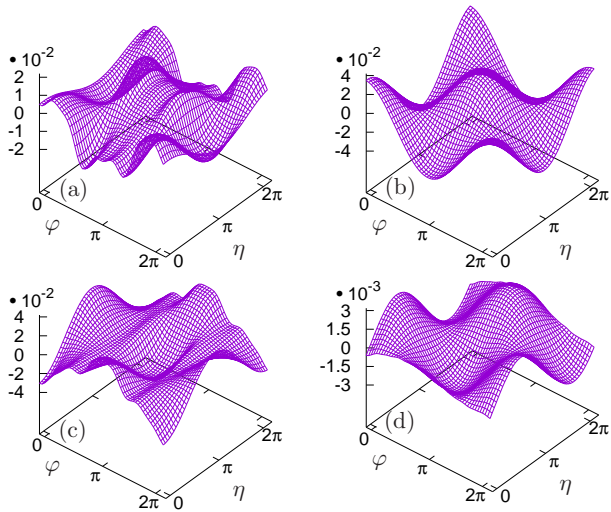}
 \vspace{0.1cm}}
\end{figure}

Here, we present the advantage of using
the IHTE for the reconstruction of the coupling functions and the iPRC. 
As an example we consider the SL oscillator 
with harmonic driving
and parameters $r=5.6$, $\mu=8$, $\alpha=0.1$, $\nu=1$ and 
$\varepsilon=0.1$ observed via variable $X_1[.]$. The coupling function is 
reconstructed by a kernel-density fit. Namely, we use a kernel 
$\mathbb{K}(x,y)=\exp[\kappa(\cos(x)+\cos(y)-2)]$ and $\kappa=200$ 
to construct $\dot{\theta}(\varphi,\eta)$. We apply a simple iterative method 
described in \cite{Kralemann_etal-13}. After $K$ iterative steps, the extracted coupling 
function $\tilde{Q}_{K,n}(\varphi,\eta) := \dot{\theta}_n(\varphi,\eta)-\tilde{\omega}$ is 
factorized into $\tilde{Z}_K(\varphi)$ and $\tilde{P}_K(\eta)$. In Fig.\ \ref{fig 3}, 
the improvement due to IHTE is evident. We used $K=30$ factorization 
steps and recover the actual coupling function with pretty high accuracy 
for different frequencies of forcing depicted in Fig.\ \ref{fig 3.1}.  

\subsection{Noisy oscillations}

\begin{figure}[h!tbp]
  \psfrag{time}[cc][cc][1.0][0]{\raisebox{0.05cm}{Time}}
  \psfrag{ylabIF}[cc][cc][1.0][0]{$\dot{\varphi}(t)$, $\dot{\theta}_{1,10}(t)$}
  \psfrag{stepn}[cc][cc][1.0][0]{\raisebox{0.0cm}{$n$}}
  \psfrag{errIF}[cc][cc][1.0][0]{\raisebox{0.0cm}{STD$^{\dot{\theta}}_n$}}
  \psfrag{Hx}[cc][cc][1.0][0]{\raisebox{0.0cm}{$\hat{H}[X_1]{\theta_{0,10}}$}}
  \psfrag{x}[cc][cc][1.0][0]{\raisebox{0.0cm}{$X_1(\theta_{0,10})$}}    
  \psfrag{c}[cc][cc][1.0][0]{\raisebox{0.0cm}{\hspace{0.0cm} (a)}}      
  \psfrag{e}[cc][cc][1.0][0]{\raisebox{0.0cm}{\hspace{0.0cm} (b)}}      
  \psfrag{h}[cc][cc][1.0][0]{(c)}      
  \psfrag{i}[cc][cc][1.0][0]{(d)}      
  \psfrag{j}[cc][cc][1.0][0]{}      
  \psfrag{k}[cc][cc][1.0][0]{(e)}      
  \includegraphics[trim = 0cm 0cm 0cm 0cm, width=0.95\columnwidth]{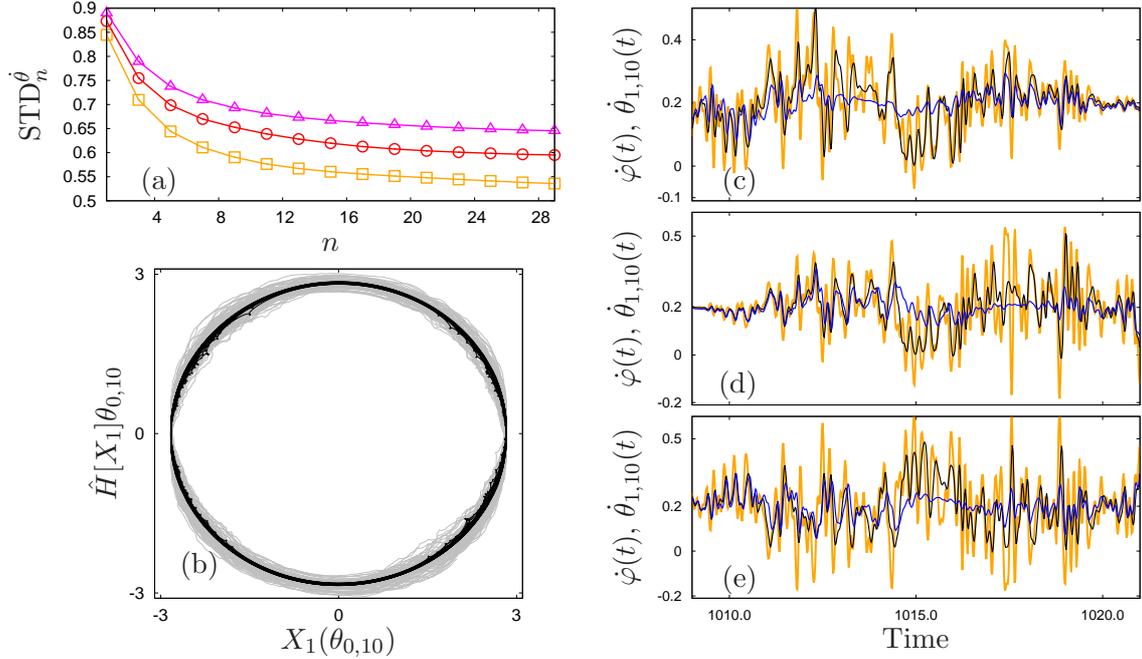}
  \caption{
  	Phase reconstruction for the SL system with $\mu=8$, $\alpha=0.1$, $\nu=1$, 
	$\varepsilon=0.2$ and $r=5.6$ and observable 
  	$X_1[.]$. Panel (a):  the errors of frequency reconstruction for 
	$\sigma= [0.06,0.08,0.1]$ (squares, circles, triangles). 
  	Panel (b): The first Hilbert embedding (grey) and the last embedding at 10th step (black). Note that 
  	the small scale loops as a result of the remaining noise influence in $X_1(t)$. Panels (c,d,e): 
  	Depicted are snapshots of the instantaneous frequency $\dot{\varphi}(t)$ (orange), and of the 
  	reconstructed frequencies $\dot{\theta}_1(t)$ (blue) 
	and $\dot{\theta}_{10}(t)$ (black) 
  	for $\sigma=[0.06,0.08,0.1]$, respectively. 
  	Note that $\dot \varphi(t)$ can be negative as an effect of noise, 
	while all reconstructions obey \eqref{eq constraints}-[I].}
 \label{fig 4}
 \end{figure}

In this section we discuss applicability of the described method to noisy signals.
We assume that the SL oscillator is driven by an external force
containing a deterministic and a stochastic (white noise) component 
\begin{equation}
\varepsilon P(t) = \varepsilon \cos(\omega r t) + \xi(t), \qquad 
\langle \xi \rangle = 0, \langle \xi(t),\xi(t')\rangle = \sigma^2 \delta(t-t')
\end{equation}
with $\mu=8$, $\alpha=0.1$, $\nu=1$, $\varepsilon=0.2$, $r=5.6$ and different noise levels 
$\sigma=[0.1,0.08,0.06]$. 
We assume a ``perfect'' observation according to $X_1$ (i.e., there is no observational noise). 
Due to the stochastic forcing, the signal's spectrum has infinite support.
In the time domain, $X(t)$ contains an infinite amount of local maxima 
and minima which will cause infinitely many but small loops 
in the embedding (see Fig.\ \ref{fig 4} (b)). Strictly speaking, we can not 
obtain phase from such a signal by calculating the length of the embedded curve, 
because the latter is a fractal curve. 

Therefore, we can not deal with the raw signal $X(t)$. Instead, as a preprocessing, 
we smooth out fast small-scale fluctuations of $X(t)$ by a SG [12,25,4] filter, 
effectively cutting the spectrum of the signal at high frequencies. 
In such a setting with a finite-width spectrum, we expect that IHTE can improve the phase 
reconstruction. 

The results in this case have to be interpreted relative to the 
smoothing parameters which are chosen in such a way that they preserve essential 
local features of the dynamics. Indeed, we observe negative 
instantaneous frequencies $\dot{\varphi}$ pointing to the need of a 
high polynomial order of smoothing (see Fig.\ \ref{fig 4} (c,d,e)). 
Also, the noise causes diffusion of phase (see Fig.\ \ref{fig 4}).

From the viewpoint of phase extraction via embeddings, 
the white noise 
forcing represents a ``worst case''. On the contrary, in all situation where coloured 
noise with a bounded spectrum is present,  we expect IHTE to be 
easier applicable. Depending on the spectral composition of noise, small-scale loops
in the embedding may be not present at all, or may be eliminated with minimal filtering.
If the noise has only relatively low-frequency component, the embedding will be
relatively smooth, and no additional processing is needed.
 
Our method is restricted to the conditions \eqref{eq constraints}. Since all of these conditions 
are not fulfilled in this example, the actual phase dynamics is only partly reconstructed, 
as can be seen also from Fig.\ \ref{fig 4} (a) where the reconstruction error decay is 
much less pronounced than in Fig.\ \ref{fig 1}. In view of this, the presented example 
can be considered as a proof of concept for IHTE of noisy signals. 
The method improves the estimation of the phase, as the examples 
of Fig.~\ref{fig 4} show, by factor up to $2$.

We add the following preprocessing to IHTE:
\begin{enumerate}
\item[1.] Given $X(t)$, apply a high-order SG-filter making the signal smooth with a large number of inflection points.
\item[2.] Next, smooth $\varphi(t)$ by the same SG-filter.
\item[3.] Proceed signal $X(t)$ with IHTE as described in Sec.\ \ref{sec: IH}.
\end{enumerate}


\section{Conclusion and open Problems}
In summary, the IHTE approach solves the problem of phase demodulation for purely 
phase modulated signals. Here, we present results for a dynamical system, where the 
amplitude dynamics is also present and linked to the dynamics of $\varphi$. We have 
demonstrated that IHTE indeed provides
a good reconstruction of the phase dynamics, if the amplitude variations are relatively 
small (see Fig.\ \ref{fig 0.1}, \ref{fig 1}, \ref{fig 2}). We show that iterations 
drastically improve the reconstruction of the phase, in comparison to the previously employed 
approach based on a single Hilbert transform 
(see Fig. \ref{fig 2}) and $Z(\varphi)$ (see Fig.\ \ref{fig 3}). However, the 
analysis of the performance of IHTE in the case of larger amplitude variations is a question 
to be discussed in the future. 

An important issue in the phase reconstruction is the protophase-to-phase transformation.
It is particularly relevant for generic observables like $X_{3}[.]$, with complex waveforms.
While handling such observables in the framework of IHTE does not state a problem,
influence of amplitude variations may depend drastically on the complexity of the waveform.
It should be stressed here, that while construction of the protophase via IHTE is almost exact,
the protophase-to-phase transformation is based on some assumption about the dynamics, which 
typically are only approximately fulfilled. This topic certainly deserves further studies.
 
Biological systems are noisy. We have given an example here, where IHTE also improves the 
reconstruction of the phase in presence of fluctuations  (see Fig.\ \ref{fig 4}). However, the very
concept of a monotonously growing phase should be reconsidered for noisy signals. 
Here we largely avoided problems by smoothing the observed signal, but in this 
approach some features of the modulation might be lost.

\section*{Acknowledgement}

Both authors thank Aneta Stefanovska and Peter McClintock for the kind invitation to contribute to this interdisciplinary work. Erik Gengel thanks the Friedrich-Ebert-Stiftung for financial support. This paper was supported
by the RSF grant 17-12-01534. The analysis in Sec. 3.3 was supported by the
Laboratory of Dynamical Systems and Applications NRU HSE, of the Russian Ministry of
Science and Higher Education (Grant No. 075-15-2019-1931).

\bibliography{pobo.bib}

\end{document}